\documentclass[dvips,12pt,a4paper]{article}
\usepackage{amsmath,amssymb}        
\usepackage{latexsym}
\usepackage{graphicx}
\usepackage{psfrag}
\usepackage{url}

\def\RR{\hbox{I\kern-.2em\hbox{R}}}
\def\ds{\displaystyle}

\title{A Non-Standard Finite Difference Scheme for MHD Boundary Layer Fluid Flow 
}

\author{Riccardo Fazio and Alessandra Jannelli  \\
Department of Mathematics, Computer Science,\\ Physical Sciences and Earth Sciences\\
University of Messina, \\
Viale F. Stagno D'Alcontres 31, 98166 Messina, Italy. \\ e-mail: rfazio@unime.it \ \ \ e-mail: ajannelli@unime.it}
\date{\today}

\begin{document}
\maketitle

\begin{abstract}
This paper deals with a non-standard finite difference scheme defined on a quasi-uniform mesh for approximate solutions of the MHD boundary layer flow of an incompressible fluid past a flat plate for a wide range of the magnetic parameter.
The obtained numerical results are compared with those available in the literature. 
We show how to improve the obtained numerical results via a mesh refinement and a Richardson extrapolation.
\end{abstract}



\section{Introduction}\label{S:Intro}
The simplest example of the application of the boundary layer theory is related to the celebrated Blasius \cite{Blasius:1908:GFK}\ problem. 
This problem describes the flow around a very thin flat plate.

The first goal of this paper is to solve numerically, with a great accuracy, the MHD boundary layer equation governing the flow of an incompressible fluid past a flat plate by a non-standard finite difference scheme defined on a quasi-uniform mesh.
Numerical methods for problems like the one considered in this paper can be classified according to the numerical treatment of the boundary condition imposed at infinity.
The oldest and simplest treatment is to replace infinity with a suitable finite value, the so-called truncated boundary.
However, being the simplest approach this has revealed within the decades some drawbacks that suggest not to apply it specially if we have to face a given problem without any clue on its solution behaviour.
Several other treatments have been proposed in literature to overcome the shortcomings  of the truncated boundary approach.
In this research area they are worth of consideration: the formulation of so-called asymptotic boundary conditions by de Hoog and Weiss 
\cite{deHoog:1980:ATB}, Lentini and Keller \cite{Lentini:BVP:1980} and Markowich \cite{Markowich:TAS:1982,Markowich:ABV:1983}; the reformulation of the given problem in a bounded domain as studied first by de Hoog and Weiss and developed more recently by Kitzhofer et al. \cite{Kitzhofer:2007:ENS}; the free boundary formulation proposed by Fazio \cite{Fazio:1992:BPF} where the unknown free boundary can be identified with a truncated boundary;  the treatment on the original domain via pseudo-spectral collocation methods, see the book by Boyd \cite{Boyd:2001:CFS} or the review by  Shen and Wang \cite{Shen:SRA:2009} for more details on this topic; and, finally, a non-standard finite difference scheme on a quasi-uniform grid defined on the original domain by Fazio and Jannelli \cite{Fazio:2014:FDS}. 
This non-standard finite difference scheme has been successively modified by Fazio and Jannelli \cite{Fazio:2017:BII}.

This study concludes by comparing the current numerical results with those given by the integral approximation method (ITM) and the non integral technique (NIT) used by Singh and Chandarki \cite{Singh:2012:NIT}.

\section{Model problem}\label{S:Model}
We consider a steady two-dimensional flow of a viscous fluid on a flat plate in the presence of a given transverse magnetic field with small electric conductivity and large transverse magnetic field. Introducing appropriate similarity variables, 
the governing equations can be reduced to the following boundary value problem (BVP) \cite{Singh:2012:NIT}
\begin{eqnarray}\label{eq:EBE}
& {\displaystyle \frac{d^3 u}{dx^3}} + u {\displaystyle \frac{d^2u}{dx^2}} + \beta \left(1 - {\displaystyle \frac{du}{dx}}\right)
= 0 \\ \nonumber
& \\
& u(0) =  {\displaystyle \frac{du}{dx}}(0) = 0 \ , \qquad \qquad {\displaystyle \frac{du}{dx}}(\infty) = 1 \ , \nonumber
\end{eqnarray}
where $\beta$ is the magnetic parameter.

\section{The finite difference scheme}\label{S:FD}
Without loss of generality we consider the class of BVPs 
\begin{eqnarray}
&& {\ds \frac{d{\bf u}}{dx}} = {\bf f} \left(x, {\bf u}\right)
\ , \quad x \in [0, \infty) \ , \nonumber \\[-1.5ex]
\label{p} \\[-1.5ex]
&& {\bf g} \left( {\bf u}(0), {\bf u} (\infty) \right) = {\bf 0}
\ ,  \nonumber
\end{eqnarray}
where $ {\bf u}(x) $ is a $ d-$dimensional vector with $ {}^{\ell} u
(x) $ for $ \ell =1, \dots , d $ as components, $ {\bf f}:[0,
\infty) \times \RR^d \rightarrow~\RR^d $, and $ {\bf g}:
\RR^d \times \RR^d \rightarrow \RR^d $.
Here, and in the following, we use Lambert's notation for the vector components \cite[pp. 1-5]{Lambert}.

In order to solve problem (\ref{p}) on the original domain we discuss first quasi-uniform grids maps from a reference finite domain and introduce on the original domain a non-standard finite difference scheme that allows us to impose the given boundary conditions exactly.
Let us consider the smooth strict monotone quasi-uniform maps $x = x(\xi)$, the so-called grid generating functions,
see Boyd \cite[pp. 325-326]{Boyd:2001:CFS} or Canuto et al. \cite[p. 96]{Canuto:2006:SMF},
\begin{equation}\label{eq:qu1}
x = -c \cdot \ln (1-\xi) \ ,
\end{equation}
and
\begin{equation}\label{eq:qu2}
x = c \frac{\xi}{1-\xi} \ ,
\end{equation}
where $ \xi \in \left[0, 1\right] $, $ x \in \left[0, \infty\right] $, and $ c > 0 $ is a control parameter.
So that, a family of uniform grids $\xi_n = n/N$ defined on interval $[0, 1]$ generates one parameter family of quasi-uniform grids $x_n = x (\xi_n)$ on the interval $[0, \infty]$.
The two maps (\ref{eq:qu1}) and (\ref{eq:qu2}) are referred as logarithmic and algebraic map, respectively. 
As far as the authors knowledge is concerned, van de Vooren and Dijkstra \cite{vandeVooren:1970:NSS} were the first to use these kind of maps. 
We notice that more than half of the intervals are in the domain with length approximately equal to $c$ and 
$x_{N-1} = c \ln N$ for (\ref{eq:qu1}),
while $ x_{N-1} \approx c N $ for (\ref{eq:qu2}).
For both maps, the equivalent mesh in $x$ is nonuniform with the
most rapid variation occurring with $c \ll x$.
The logarithmic map (\ref{eq:qu1}) gives slightly better resolution near $x = 0$ than the
algebraic map (\ref{eq:qu2}), while the algebraic map gives much better resolution than the
logarithmic map as $x \rightarrow \infty$. 
In fact, it is easily verified that
\[
-c \cdot \ln (1-\xi) < c \frac{\xi}{1-\xi} \ ,
\]
for all $\xi$, but $\xi = 0$.

The problem under consideration can be discretized by introducing a uniform grid $ \xi_n $ of $N+1$ nodes in $ \left[0, 1\right] $ 
with $\xi_0 = 0$ and $ \xi_{n+1} = \xi_n + h $ with $ h = 1/N $, so that $ x_n $ is a quasi-uniform grid in $ \left[0, \infty\right] $. 
The last interval in (\ref{eq:qu1}) and (\ref{eq:qu2}), namely $ \left[x_{N-1}, x_N\right] $, is infinite but the point $ x_{N-1/2} $ is finite, because the non integer nodes are defined by 
\[
x_{n+\alpha} = x\left(\xi=\frac{n+\alpha}{N}\right) \ ,
\]
with $ n \in \{0, 1, \dots, N-1\} $ and $ 0 < \alpha < 1 $. 
These maps allow us to describe the infinite domain by a finite number of intervals.
The last node of such grid is placed on infinity so right boundary conditions
are taken into account correctly.

We approximate the values of the scalar variable $u(x)$ and its derivative at mid-points of the grid $x_{n+1/2}$, for $n=0,\cdots,N-1$, 
using non-standard difference discretizations
\begin{eqnarray}
u_{n+1/2} &\approx& \frac{x_{n+3/4}-x_{n+1/2}}{x_{n+3/4}-x_{n+1/4}} u_n + \frac{x_{n+1/2}-x_{n+1/4}}{x_{n+3/4}-x_{n+1/4}} u_{n+1} \ , \nonumber \\
[-1.5ex] \label{eq:udu} \\ [-1.5ex]
\left. \frac{du}{dx}\right|_{n+1/2} &\approx& \frac{u_{n+1}-u_n}{2\left(x_{n+3/4} - x_{n+1/4}\right)} \nonumber \ .
\end{eqnarray}
We emphasize that the key advantage of our non-standard finite difference formulation is to overcome the 
difficulty of the numerical treatment of the boundary conditions at the infinity. 
In fact, the formulae (\ref{eq:udu}) use the value $ u_N = u(\infty) $, but not $ x_N = \infty $
and then, the boundary conditions at infinity are taken into account in a natural way.

For the class of BVPs (\ref{p}), a non-standard finite difference scheme on a quasi-uniform grid can be defined by using the approximations 
given by (\ref{eq:udu}) above, and it can be written as follows
\begin{eqnarray}
& {\bf U}_{n+1} - {\bf U}_{n} - a_{n+1/2} {\bf f} \left( x_{n+1/2}, b_{n+1/2}{\bf U}_{n+1} + c_{n+1/2}{\bf U}_{n} \right) = {\bf 0}
\ , \nonumber\\
& \mbox{for} \quad n=0, 1, \dots , N-1
\label{boxs} \\ 
& {\bf g} \left( {\bf U}_0,{\bf U}_N \right) = {\bf 0} \ ,  \nonumber
\end{eqnarray}
where 
\begin{eqnarray}\label{eq:abc}
a_{n+1/2} &=& 2\left(x_{n+3/4} - x_{n+1/4}\right) \ , \nonumber \\
b_{n+1/2} &=&  \frac{x_{n+1/2} - x_{n+1/4}}{x_{n+3/4}-x_{n+1/4}} \ , \\
c_{n+1/2} &=&  \frac{x_{n+3/4} - x_{n+1/2}}{x_{n+3/4}-x_{n+1/4}} \nonumber \ , 
\end{eqnarray}
for $n=0, 1, \dots , N-1$.
The finite difference formulation (\ref{boxs}) has order of accuracy $O(N^{-2})$. It is evident that (\ref{boxs}) is
 a nonlinear system of $ d \; (N+1)$ equations in the $ d \; (N+1)$ unknowns $ {\bf U} = ({\bf U}_0,{\bf U}_1, \dots , {\bf U}_N)^T $. 
For the solution of (\ref{boxs}) we can apply the classical Newton's method along with the simple termination criterion
\begin{equation}\label{eq:Tcriterion}
{\ds \frac{1}{d(N+1)} \sum_{\ell =1}^{d} \sum_{n=0}^{N}
|\Delta {}^\ell U_{n}| \leq {\rm TOL}} \ ,
\end{equation}
where $ \Delta {}^ \ell U_{n} $, for $ n = 0,1, \dots, N $ and $ \ell = 1, 2, \dots , d $, is the difference between two successive iterate components and $ {\rm TOL} $ is a fixed tolerance. 

\section{Numerical results and comparison}\label{S:Results}
In this Section, we present the numerical results obtained by solving the mathematical model (\ref{eq:EBE}) 
using the non-standard finite difference scheme (\ref{boxs}) on the quasi-uniform grid defined by the logarithmic map (\ref{eq:qu1}) with control parameter $c=2$. 
Now, let us rewrite the model (\ref{eq:EBE}) as a first order system as follows
\begin{eqnarray}
&& \frac{d {}^1u}{dx}= {}^2u,  \nonumber \\
&& \frac{d {}^2u}{dx}= {}^3u,  \qquad x \in ( 0, \infty ) \label{sys:EBE} \\
&& \frac{d {}^3u}{dx}= - {}^1u {}^3u  - \beta (1 - {}^2u ),\nonumber
\end{eqnarray}
with
\begin{eqnarray*}
\left.
\begin{array}{ll}
{}^1u(0)= {}^2u(0)=0 \ , \qquad {}^2u(\infty)=1 \ , \label{bv:EBE} \\
\end{array}
\right.
\end{eqnarray*}
or, in an equivalent form, 
\begin{eqnarray*}
&& \textbf{u}=({}^1u,{}^2u,,{}^3u)^T , \nonumber \\
&& \textbf{f}( x,\textbf{u})= \left( {}^2u, {}^3u, - {}^1u {}^3u  - \beta (1 - {}^2u ) \right)^T, \\
&& \textbf{g}(\textbf{u}(0),\textbf{u}(\infty))=({}^1u(0),{}^2u(0),{}^2u(\infty)-1)^T \ , \nonumber
\end{eqnarray*}
where $\textbf{u}(x)$ is a three-dimensional vector with components ${}^\ell \emph{u}(x)$ for $\ell=1,2,3$, and
$\textbf{f}:[0,\infty) \times  \mathbb{R}^d\to  \mathbb{R}^d$ and $\textbf{g}: \mathbb{R}^d \times  \mathbb{R}^d\to \mathbb{R}^d$,
with $d=3$. We set as first guess for the Newton's iteration, and for the whole range of $\beta$, with $\beta=0,0.2,\cdots,2$, 
the following initial data
\begin{eqnarray}\label{EBE_IC}
&& {}^1u(x)=0.5 \ x  \ , \qquad {}^2u(x)= 1  \ ,  \qquad  {}^3u(x) = \exp(- x). \nonumber
\end{eqnarray}
Moreover, for all tests we consider a fixed tolerance ${\rm TOL} = 10^{-8}$ and $N=1000$.

\begin{figure}[!ht]
\begin{center}
\hfil
\includegraphics[width=.85\textwidth]{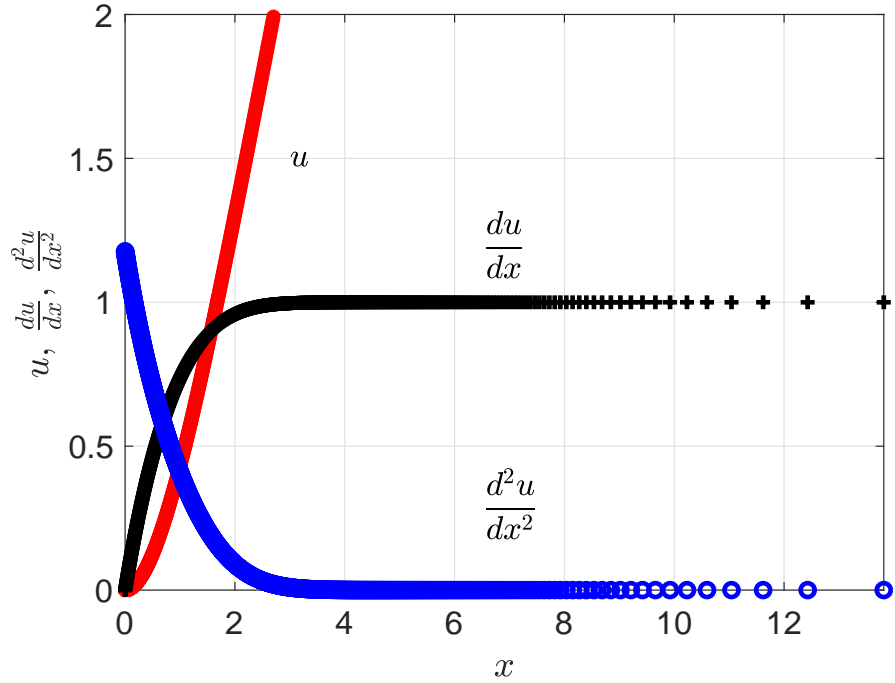}
\hfil 
\end{center}
\caption{Numerical solution for the problem (\ref{eq:EBE}) for $\beta=1.2$.}
\label{solution}
\end{figure}

In Figure \ref{solution}, we report the numerical solution obtained for $\beta=1.2$.
The recovered value of the second order derivative of the solution at the origin is 
${\displaystyle \frac{d^2u}{dx^2}}(0)=1.177226684282633$, obtained in $6$ iterations.

The table \ref{tab_comp} lists the obtained numerical results. 
For the sake of brevity we have chosen to report only the values of the wall shear stress, that is the second derivative value at the origin.
Within the same table we can compare our results with those reported by Singh and Chandarki \cite{Singh:2012:NIT}.
The problem with $\beta = 0 $ corresponding to the Blasius problem.
\begin{table}
\centering
\begin{tabular}{cccc}
\hline
$ \beta $ & DTM \cite{Singh:2012:NIT} & NIT \cite{Singh:2012:NIT} & FD (this study)
\\ 
\hline
0.0 & 0.46910 & 0.46920 & 0.4695998 \\
0.2 & 0.66343 & 0.64819 & 0.6389912 \\
0.4 & 0.80009 & 0.78749 & 0.7749667 \\
0.6 & 0.91659 & 0.90562 & 0.8917423 \\
0.8 & 1.01988 & 1.01002 & 0.9956201 \\
1.0 & 1.11362 & 1.10460 & 1.0900651 \\
1.2 & 1.20006 & 1.19170 & 1.1772267 \\
1.4 & 1.28068 & 1.27285 & 1.2585472 \\
1.6 & 1.35652 & 1.34913 & 1.3350501 \\
1.8 & 1.42834 & 1.42132 & 1.4074922 \\
2.0 & 1.49671 & 1.49002 & 1.4764520 \\
\hline
\end{tabular}
\caption{Numerical results, related to $\frac{d^2u}{dx^2}(0)$, and comparison.}
\label{tab_comp}
\end{table}

We improve the accuracy of the computed solution through
subsequent refinements of the computational domain by using the Richardson's extrapolation.
On the computational domain of the problem, we build a quasi-uniform grid with a mesh-points number equal 
to $ N_0 $ and proceed with subsequent grid refinements by constructing meshes with grid-point numbers $ N_g $ for $ g = 1, 2, \cdots$,
where $ N_ {g + 1} = r N_ {g} $ with refinement factor $r=2$.
On each grid, the numerical solution $ U_g $, $ g = 0,1, \cdots, G $ is computed using the non-standard finite difference method.
In order to reduce the calculations, we adopt a continuation strategy, in fact we use the final solution  $ U_g $
obtained on the grid $ g $ as initial guess for calculating the solution $ U_{g+1}$ on the grid $g+1$.
where the new grid values are approximated by linear interpolations.
We define the level of the Richardson's extrapolation by the index $k$ and, 
the two numerical solutions related to the grids $g$ and $g+1$ at the extrapolated level $k$ by $ U_{g,k} $ and $ U_{g+1,k}$. 
We use the following formula to calculate a more accurate approximation 
\begin{equation} \label{eq:Rextra}
U_{g+1,k+1} = U_{g+1,k} + \frac{U_{g+1,k}-U_{g,k}}{2^{p_k}-1} \qquad k=0,1,\cdots,G-1 \ .
\end{equation}
In table \ref{Tab_extr}, we report the extrapolated values with $N=100,200,400$ grid points for $\beta=1$. 
The last extrapolated value is  $^3U_{2,2} = 1.090064908$ and can be considered as our benchmark value for 
${\displaystyle \frac{d^2u}{dx^2}}(0)$. 
We can conclude that the reported extrapolated value is correct up to $9$ decimal places.

\begin{table}[!h]
\centering
\begin{tabular}{rrrr}
\hline
\hline
$N_g$&$^3U_{g,0}$&$^3U_{g,1}$&$^3U_{g,2}$\\ 
\hline\hline
 $100$  &  $1.090081494$  & ---  \\
 $200$  &  $1.090069055$  & $1.090064908$ & --- \\
 $400$  &  $1.090065945$  & $1.090064908$ & $1.090064908$ \\
\hline
\hline
\end{tabular}
\caption{Extrapolated values at origin $x=0$ for ${\displaystyle \frac{d^2u}{dx^2}}(0)$ with $\beta=1$ }
\label{Tab_extr}
\end{table}

\section{Concluding Remarks}\label{S:Remarks}
In this paper the problem \ref{eq:EBE}, that describes the MHD boundary layer flow of  an incompressible fluid past a flat plate, is solved by non-standard finite difference method 
on quasi-uniform grid for the different magnetic parameters $\beta$. 
The values of the second order derivative of the solution at the origin for different values of parameter $\beta$
are reported in the Table \ref{tab_comp}. In order to verify the accuracy
of the proposed method, the results are compared with those by Singh and Chandarki \cite{Singh:2012:NIT}.
The recovered values are accurate.

\noindent {\bf Acknowledgments.} {The research of this work was supported, in part, 
by the University of Messina and by the GNCS of INDAM.}

\end{document}